\documentclass[oneside, a4paper,11pt,reqno]{amsart}
\usepackage{color}
\usepackage{amssymb,amsmath,amsthm}
\usepackage{graphics,graphicx}
\usepackage[T1]{fontenc}

\newtheorem{theorem}{Theorem}[section]
\newtheorem{lemma}[theorem]{Lemma}
\newtheorem{e-proposition}[theorem]{Proposition}

\newtheorem{e-definition}[theorem]{Definition\rm}

\def\A{\mathcal{A}}
\def\B{\mathcal{B}}
\def\C{\mathcal{C}}

\def\S{\mathcal{S}}

\def\E{\mathcal{E}}
\def\R{\mathcal{R}}
\def\F{\mathcal{F}}
\def\T{\mathcal{T}}

\def\class{\textrm{class}}
\def\branch{\textrm{Branch}}

\def\PP{\mathbb{P}}

\def\CC{\mathbb{C}}

\title{Degree and class of caustics by reflection for a generic source}
\date\today

\author{Alfrederic Josse}
\address{Universit\'e de Brest,
LMBA, UMR CNRS 6205, 29238 Brest cedex, France}
\email{alfrederic.josse@univ-brest.fr}

\author{Fran\c{c}oise P\`ene}
\address{Universit\'e de Brest,
LMBA, UMR CNRS 6205, 29238 Brest cedex, France}
\email{francoise.pene@univ-brest.fr}

\keywords{caustic, degree, class, birationality\\
Fran\c{c}oise P\`ene is supported by the french ANR project GEODE 
(ANR-10-JCJC-0108).}

\begin{document}

\begin{abstract}
Given any irreducible algebraic (mirror) curve $\C\subseteq \PP^2:=\PP^2(\CC)$ and any
(light position) $ S \in\PP^2$, the caustic by reflection $\Sigma_{S}(\C)$ of $\C $ from $ S $
is the Zariski closure of the envelope of the reflected lines got from the lines coming from $S$
after reflection on $\C $. In \cite{fredsoaz1,fredsoaz2}, we established formulas
for the degree  and class (with multiplicity) of $\Sigma_{S}(\C)$
for any $\C $ and any $ S $. In this paper, we prove the birationality
of the caustic map for a generic $S$ in $\PP^2$. Moreover, we give 
simple formulas for the degree and class (without multiplicity)
of $\Sigma_{S}(\C)$ for any $\C $ and for a generic $S$ in $\PP^2 $.




%

\end{abstract}
\maketitle

\section{Introduction}
\label{}
We are interested in the study of caustics by
reflection in the projective complex plane
$\PP^2$.
Given an irreducible algebraic curve $\C =V (F)\subset \PP^2$ of degree $ d\ge 2 $ and
given $S=[ x_0:y_0:z_0 ]\in\PP^2 $, the {\bf caustic by reflection} $\Sigma_{S}(\C)$ of $\C $ from $S $ is
the Zariski closure of the envelope of the reflected lines on $\C $
of the lines coming from $S $.

For $ m\in\C $, the reflected line $\R_{m,S,\C}$ is defined as the
orthogonal symmetric of the (incident) line $(m\, S)$ with respect to
the tangent line to $\C$ at $m$.
In \cite{fredsoaz1,fredsoaz2}, we detail the construction of the reflected lines and
we define two rational maps $ \rho _{F,{S}}$ and $\Phi_{F,{S}}$ from $\PP^2 $
into itself satisfying the following property:
For a generic $ m $ in $\C$, $\rho_{F,{S}}(m)$ corresponds to an equation of the
reflected line $\R_{m,S,\C}$  and this line is tangent to $\Phi_{F,{S}}(\C)$
at $\Phi_{F,{S}}(m)$.
Hence the caustic $\Sigma_{S}(\C)$ is the Zariski closure of
$\Phi_{F,{S}}(\C)$ and $\Phi_{F,{S}}$ is called the {\bf caustic map} of $\C$ 
from $S$. Observe that the Zariski closure of $\rho_{F,{S}}(\C)$
is then the dual curve of the caustic $\Sigma_{S}(\C)$.
In \cite{fredsoaz1,fredsoaz2}, we used this approach
to establish precise formulas for the degree anf class (both with multiplicity) of $\Sigma_{S}(\C)$
for any $\C $ and any $S $. The degree with multiplicity of $\Sigma_{S}(\C)$ means its degree multiplied
by the degree of the rational map $\Phi_{F, {S} }$
restricted to $\C$ . The class with multiplicity of $\Sigma_{\S}(\C)$ means its class
multiplied by the degree of the rational map $\rho_{F,{S}}$ restricted to $\C$.
Our formulas complete the formula obtained by Chasles in \cite{Chasles} for the class of a caustic by reflection
(for a generic $\C$ and a generic $S$). Let us indicate that, in \cite{Brocard-Lemoyne}, Brocard and Lemoyne gave, without
any proof, formulas for the degree and class of caustics by reflection (for a Pl\"ucker curve $\C$ and for $S$ not at infinity).
It seems that their formulas come from an incorrect composition of formulas by Salmon and Cayley \cite{Salmon-Cayley} for some
caracteristic invariants of
pedal and evolute curves (using the representation of caustics by reflection
given by Quetelet and Dandelin). This is discussed in \cite{fredsoaz2}.
Let us also mention the work of Catanese and Trifogli on focal loci,
which generalize evolutes to higher dimension \cite{Trifogli,CataneseTrifogli}.

The question of the birationality of the rational maps $\rho_{F,{S}}$ and $\Phi_{F,{S}}$
on $\C$
is not evident even if $ S $ is not at infinity.
Indeed, according to results of Quetelet and Dandelin \cite{Quetelet,Dandelin}, when $ S $ is not at infinity,
the caustic $\Sigma_{S}(\C)$ is the evolute of the $S$-centered homothety (with
ratio 2) of the pedal  of $\C$ from $S$
(i.e. the evolute of the orthotomic of $\C$ with respect to $S$).
But we just know that the evolute map is birational for a generic algebraic curve
(see \cite{Fantechi} by Fantechi).

In this note, we prove the birationality on $\C$ of the maps $ \rho _{F,{S}}$ and $\Phi_{F,{S}}$
for any irreducible algebraic curve $\C \subset \PP^2$ of degree $d\ge 2$ and for a generic $S$ in $\PP^2$.
This result enables us to establish simple
formulas for the degree and class of caustics
by reflection valid for any irreducible algebraic curves $\C\subset \PP^2 $ of degree $d\ge 2$ and for a generic $ S $ in
$\PP^2$. In this study,  the  
cyclic points $I=[1:i:0]$ and $J=[1:-i:0]$ play a particular role. We 
will also use the canonical projection $\pi:\CC^3\setminus\{0\}\rightarrow \PP^2$.
\section{Birationality}
\begin{theorem}\label{thm}
Let $\C =V (F)\subset \PP^2$ be any irreducible algebraic curve of degree $ d\ge 2 $.
For a generic $S\in\PP^2 $, the maps $ \rho_{F,{S}}$ and $\Phi_{F,{S}}$ are birational on $\mathcal C $.
\end{theorem}
Before going into the proof of our Theorem, let us
introduce some notations and recall some facts (see \cite{fredsoaz1}).
For any line $ \mathcal D=V (ax+by+cz)\in\PP^2 $ such that $ a^2+b^2\ne 0 $,
we define the orthogonal symmetric with respect to $\mathcal D $ as the rational map
$\sigma_{\mathcal D }:\PP^2\rightarrow \PP^2$ (which is an involution) given by
$$\sigma_{\mathcal D }[x: y: z]=\pi\left((a^2+b^2)\cdot (x, y, z)+(ax+by+cz)\cdot (a, b, 0)\right) .$$
Let $\C =V (F)\subset\PP^2$ be an irreducible algebraic curve of degree $ d\ge 2$ and let $ S\in\PP^2 \setminus\{I, J\}$
be a (light) position.
We define $C_0:=\C \setminus V (F_x^2+F_y^2) $. Observe that this
set corresponds to the complement in $\C$ of the cyclic apparent contour of $\C $
(the cyclic apparent contour of $\C$ from the cyclic points).
We recall that the reflected line $\mathcal R_{m, S,\mathcal C}$ at $ m\in \C_0\setminus \{S\}$ is the line
$(m\, \sigma_{\mathcal T_m\C}(S)) $, where $\mathcal T_m\C $ is the tangent to $\C $ at $ m $.
For any $ m\in\C_0 $, we define the normal line $\mathcal N_m\C $ to $\C $ at $ m $ as the line
containing $ m $ and $[F_x(m):F_y(m):0]$.

For any $m\in\C_0$, we consider the set $ K_m $ of points $ S \in\PP^2$ such that there exists
$ m'\in\C_0\setminus \{m\}$ satisfying $ \rho _{F,{S}}(m')=  \rho _{F, {S} }(m) \ne 0$.
Observe that the set $\mathcal A $ of $ S\in\PP^2 $ such that $ \rho _{F, {S} }$ is not birational can be written
$\mathcal A=\bigcup_{E\subset \C_0\,:\,\# E <\infty}\bigcap_{m\in\C_0\setminus E}K_m.$
To prove that $ \rho _{F,{S}} $ is
birational for a generic $S$ in $\PP^2$, we prove that $\A$ is contained in a subvariety of
codimension at least 1 in $\PP^2$.
Our proof is based on the following
lemma.
\begin{lemma}\label{lem}
For any $ m\in \C_0 $, the set $ K_m $ is contained in a (possibly non irreducible)
algebraic curve $\bar K_m$ of degree at most $2d^2+2$.
\end{lemma}
\noindent{Proof.\/}
Let us consider any $ m\in \C_0 $.
Let $ S\in K_m $ and $ m'\in\C_0\setminus \{m\}$ satisfying $ \rho _{F, {S} }(m')= 
\rho _{F, {S} }(m) \ne 0$.
Then
$(m\, m')=\mathcal R_{m, S,\C}= \mathcal R_{m', S,\C}$
and so $ S $ is in
$ \mathcal A_{m, m'}:=
\sigma_{\mathcal T_m\C}((m\, m'))\cap
 \sigma_{\mathcal T_{m'}\C}((m\, m')).$
Observe that, if
$ \sigma_{\mathcal T_m\C}((m\, m'))
=\sigma_{\mathcal T_{m'}\C}((m\, m')),
$
then these lines are $(m\,m')$ and so $(m\,m') $ is stable by
$ \sigma_{\mathcal T_m\C} $ and by
$  \sigma_{\mathcal T_{m'}\C} $.
But $\mathcal T_m\C$ and $\mathcal N_m\C$
are the only lines containing $ m $  which are stable by $ \sigma_{\mathcal T_m\C} $.
Therefore,  
$\sigma_{\mathcal T_m\C}((m\, m'))
=\sigma_{\mathcal T_{m'}\C}((m\, m')),
$
implies that
$(m\, m')\in\{\mathcal T_m\C,\mathcal N_m\C\}\cap \{\mathcal T_{m'}\C,\mathcal N_{m'}\C\}.$
If
$\{\mathcal T_m\C,\mathcal N_m\C\}\cap \{\mathcal T_{m'}\C,\mathcal N_{m'}\C\}=\emptyset$,
then $ S $ is the only point of $ \mathcal A_{m,m'} $, so $S$ is equal to
$$ \tau_m(m'):=\pi\left((m\wedge \sigma_{\mathcal T_m\C}(m'))
\wedge  (m'\wedge \sigma_{\mathcal T_{m'}\C}(m))\right).$$
Notice that $\tau_m$ is a rational map with coordinates of degree $2d$.
We obtain that $ S $ belongs to the Zariski closure of
$\tau_m (\C) $, which (according to \cite[Proposition 4.4]{Fulton}) 
is contained in an algebraic curve of degree at most
$\C\cdot\tau_m^*(H)\le 2d^2$
(where $H$ is the hyperplane class in $\PP^2$).
Otherwise, $S\in\mathcal A_{m,m'}=(m\, m')\in \{\mathcal T_m\C,\mathcal N_m\C\} $.
Finally, we have
$ K_m\subseteq \bar K_m:=\overline{\tau_m (\C)}\cup \mathcal T_m\C\cup
\mathcal N_m \C $
which is an algebraic curve of degree at most $2d^2+2 $ (use
for example the fundamental lemma of \cite{fredsoaz1} applied
with $\tau_m$).
\qed

\noindent{\it Proof of Theorem \ref{thm}.\/}
Let us prove that $ \rho _{F,{S}} $ is birational on $\C $ for a generic $ S $ in $\PP^2$. The birationality of $ \Phi_{F,{S}} $ will follow.
Indeed, for a generic $ S $ in $\PP^2 $, the caustic $\Sigma_S (\C)$ is a curve (see for
example \cite{fredsoaz2}). Therefore, for generic $ m, m'\in\C $,
$ \Phi_{F,{S}}(m)= \Phi_{F,{S}}(m') $ implies that $ \rho _{F,{S}}(m)=
\rho _{F,{S}}(m') $.
With the notations of Lemma \ref{lem}, we define
$\mathcal A':=\bigcup_{E\subset \C_0\,:\,\# E <\infty}\bigcap_{m\in\C_0\setminus E}\bar K_m.$
We prove that the set 
$\F:=\left\{\bigcap_{m\in\C_0\setminus E}\bar K_m,\ E\subset \C_0,\ \# E <\infty\right\}$ 
is inductive for the inclusion.
Let
$\left (\F_j:=\bigcap_{m\in\C_0\setminus E_j}\bar K_m\right)_{j\ge 1}$
be an increasing sequence of sets belonging to $\F$.
Let us show that the union $Z$ of these  sets is also in $\F$.
First
$ Z\subseteq\bigcap_{m\in\C_0\setminus \bigcup_{i\ge 1}E_i}\bar K_m\subseteq
\bar K_{m_0}$
for some fixed
$m_0\in\C_0\setminus \bigcup_{i\ge 1}E_i$.
Now $\bar K_{m_0}$ is the union of a finite number of irreducible algebraic curves
$C_1,\dots, C_p$.
Let $i\in\{1,\dots, p\}$ and let $d_i$ be the degree of $C_i$.
If $C_i\subseteq Z$, then there exists
$N_i\ge 1$ such that $C_i\subseteq \F_{N_i}$.
Assume now that $C_i\not\subseteq Z$.
Then $(C_i\cap \F_j)_{j\ge 1}$ is an increasing sequence of
finite sets containing at most $d_i\times (2d^2+2)$ points.
Therefore, there exists
$N_i\ge 1$ such that $(C_i\cap Z)\subseteq \F_{N_i}$.
We conclude that $Z=\F_{\max(N_1,\dots, N_p)}$ and so $Z$ is in $ \F $. So $\F$ is inductive.

{}From the Zorn lemma, either
$ \F $ is empty or it
admits a maximal element (for the inclusion).
If it is empty, then $\A=\A'=\emptyset$.
If it is not empty and if $\F_0:=\bigcap_{m\in\C_0\setminus E_0}\bar K_m$
(with $E_0\subset \C_0$ and $\#E_0 <\infty$) is a maximal element of $ \F $,
then $\A'=\F_0$.
Indeed, $\A'$ contains $\F_0$ by definition of $\A'$.
Conversely, let $S\in \A'$, there exists
$E\subset \C_0$ such that $\#E <\infty$ and such that
$S\in \bigcap_{m\in\C_0\setminus E}\bar K_m$. 
Hence
$ S\in
 \bigcap_{m\in\C_0\setminus (E\cup E_0)}\bar K_m$.
Since we also have
$ \bigcap_{m\in\C_0\setminus E_0}\bar K_m\subseteq  
 \bigcap_{m\in\C_0\setminus (E\cup E_0)}\bar K_m$,
we conclude that
$S\in \bigcap_{m\in\C_0\setminus E_0}\bar K_m$. 
Therefore, in any case, $\A$ is contained in an algebraic curve, this gives the
$S$-genericity of the birationality of $ \rho _{F,{S}} $ and so the statement
of Theorem \ref{thm}.
\qed
\section{Light generic formulas for the degree and the class of caustics}
Let $\C =V (F)\subset \PP^2$ be any irreducible algebraic curve of degree $ d\ge 2 $.
We call {\bf isotropic tangent} to $\C$ any tangent to $\C$ containing $I$ or $J$.
Before stating our formulas, let us
introduce some notations.

For any $P\in\PP^2$, we write $\mu_P(\C)$ for the multiplicity of $\C $ at $ P $.
We recall that $\mu_P(\C)=1$ means that $ P $ is a non singular point of $\C $.
For any $ P\in\C $, we write $\branch_P (\C)$ for the set of branches of
$\C$ at $ P $.
Let us write $\E_{\C}$ for the set of couples $(P,\B) $ with $ P\in\C $ and with $\B
\in\branch_P (\C) $. For any $ (P,\B)\in \E_{\C}$, we write
$\T_P\B$ for the tangent line to $\B $ at $ P $ and $e_{\B}$ for the multiplicity
of $\B $. We recall that 
$\sum_{\B\in\branch_P(\C)}e_{\B}=\mu_{P}(\C)$.
For any $ (P,\B)\in \E_{\C}$ and any algebraic curve $\C'$, we denote by $i_P(\C,\C')$
(resp. $i_P(\B,\C')$) the intersection number of $\C $ (resp. $\B$) with $\C'$ at $ P $.
We recall that the contact number $\Omega_{m_1}(\C,\C')$ of $\C $ with $\C'$ is given by
$ \Omega_{m_1}(\C,\C')=i_P (\C,\C')- \mu_P(\C) \mu_P(\C')$.
The line at infinity of $\PP^2$ is written $\ell_\infty$.
Combining Theorem \ref{thm} with the main results
of \cite{fredsoaz1,fredsoaz2}, we obtain:
\begin{e-proposition}
Let $\C =V (F)\subseteq \PP^2$ be any irreducible algebraic curve of degree $ d\ge 2 $ and of
class $ d^\vee $.
For a generic $\S\in\PP^2 $, we have
$$ \deg( \Sigma_{S}(\C) )=3d+f_0-t_I-t_J\ \ 
\mbox{and}\ \  \class(\Sigma_{\S}(\C))=2d^\vee+d-g-\mu_I(\C)-\mu_J(\C),$$
where
$  g$ is the contact number of $\mathcal C$ with $\ell_\infty$, i.e.
$ g:=\sum_{m_1\in\C\cap\ell_\infty}\Omega_{m_1}(\mathcal C,\ell_\infty), $
where $ f_0 $ is the number of ``inflectional branches'' of $\C$ not tangent to the  line 
at infinity, i.e.
$$ f_0:=\sum_{(P,\B)\in \E_{\C}\,:\, i_{P}(\B,\T_P\B)>2e_{\B},\,\T_P\B\ne\ell_\infty }( i_{P}(\B,\T_P\B)-2e_{\B} ).$$
and where $ t_P$ is the number of branches of $\C $
tangent at $ P$ to the line at infinity:
$t_P=\sum_{\B\in Branch_P(\C)\,:\,\mathcal T_P\mathcal B=\ell_\infty} e_{\mathcal B}.$
\end{e-proposition}
\noindent{\it Proof.\/}
According to Theorem \ref{thm}, for a generic $ S $ in $\PP^2$, the degree
and class (with multiplicity) of $ \Sigma_{S}(\C) $ are equal to its degree
and class.
For the degree formula,  we use Theorem 20 of \cite{fredsoaz1}.
For the class formula, we use Theorem 2 of \cite{fredsoaz2}.
For a generic $ S\in\PP^2 $ ($S\in\PP^2\setminus (\C \cup \ell_\infty)$
not contained in an isotropic tangent to $\C $), we have
$f=\mu_I(\C)+\mu_J(\C)$, $f'=0$, $g'=0$ and $q'=0$
(with the notations of Theorem 2 in \cite{fredsoaz2}).
\qed
%
%
%




\section*{Acknowledgements}
We thank Fabrizio Catanese for discussions having motivated the redaction of
this note on source generic results for caustics by reflection.

\end{document}